\documentclass{amsart}
\usepackage{amssymb}  
\usepackage[all]{xy}  
  
\newcommand\bS{{}^b\kern-2ptS}  
\newcommand\PsS[2]{\Psi_{\rm inv}^{#1}(#2)}  
\newcommand\Psm[2]{\psi_{\rm inv}^{#1}(#2)}  
\newcommand\Ps[2]{\Psi^{#1}(#2)}  
\newcommand\End{\operatorname{End}}  
\newcommand\Hom{\operatorname{Hom}}  
\newcommand\alg{\operatorname{alg}}  
\newcommand\In{\operatorname{In}}  
\newcommand\GR{\mathcal{G}}  
\newcommand\pSS[1]{{}^p\kern-1pt\SS{#1}}  
\newcommand\lra{\longrightarrow}  
\newcommand\ind{\operatorname{ind}}   
\newcommand\mfk{\mathfrak}  
  
\newcommand\datver[1]{\def\datverp%  
 {\par\boxed{\boxed{\text{Version: #1; Run: \today}}}}}  
\datver{1.2; Revised: August 7, 2001}

\newcommand\CC{\mathbb C}

\newcommand\RR{\mathbb R}  
\newcommand\ZZ{\mathbb Z}  
\renewcommand\SS{\mathbb S}  
  
\newcommand\lgg{\mathfrak g}  
\newcommand\pa{\partial}

\newcommand\CI{{\mathcal C}^{\infty}}  
\newcommand\CIc{{\mathcal C}^{\infty}_{\text{c}}}  
\newcommand\ie{{\em i.e.,} }  
\newtheorem{theorem}{Theorem}  
  
\newtheorem{corollary}{Corollary}  
\newtheorem{lemma}{Lemma}  
  
\theoremstyle{definition}

\theoremstyle{remark}

\begin{document}  
\title[Index for families] {An index theorem for gauge-invariant families: 
The case of solvable groups}  
  
\author{Victor Nistor} \address{Current address:\  
  Department of Mathematics,  
  Pennsylvania State University} \email{nistor@math.psu.edu}  
\thanks{Partially supported by the NSF Young Investigator Award  
  DMS-9457859, a Sloan fellowship, and NSF Grant  
DMS-9971951. {\bf  
    http:{\scriptsize//}www.math.psu.edu{\scriptsize/}nistor{\scriptsize/}}.  
%Key words: Elliptic operator, group action, index formula, $K$-theory,  
%cyclic homology, non-commutative geometry  
%eta invariant, residues. AMS Subject Classification 1986: 46L80,58G20  
}  
\dedicatory\datverp  
%\dedicatory{\today}  

\begin{abstract} 
We define the gauge-equivariant index of a family of elliptic
operators invariant with respect to the free action of a family $\GR
\to B$ of Lie groups (these families are called ``gauge-invariant
families'' in what follows).  If the fibers of $\GR \to B$ are
simply-connected and solvable, we compute the Chern character of the
gauge-equivariant index, the result being given by an Atiyah-Singer
type formula that incorporates also topological information about the
bundle $\GR \to B$. The algebras of invariant pseudodifferential
operators that we study, $\Psm {\infty}Y$ and $\PsS {\infty}Y$, are
generalizations of ``parameter dependent'' algebras of
pseudodifferential operators (with parameter in $\RR^q$), so our
results provide also an index theorem for elliptic, parameter
dependent pseudodifferential operators.  We apply these results to
study Fredholm boundary conditions on a simplex.
\end{abstract}  
\maketitle  
  
\tableofcontents  
  
\section*{Introduction}

Families of operators invariant with respect to the action of a family
of Lie groups (referred to as ``gauge-invariant families,'' in what
follows) have been considered in several papers in the literature.
They appear, for example, in the work of Bismut on analytic torsion
\cite{BismutSURV}, in the work of Connes on abstract index theory on
groupoids \cite{ConnesNCG}, and in the work of Mazzeo and Melrose on
analysis on non-compact manifolds \cite{Mazzeo, MazzeoMelrose,
MelroseCongress}. Gauge-invariant operators are also a basic example
of the operators on groupoids introduced in \cite{NWX}.

Gauge-invariant families are likely to form the main building block in
the analysis of operators on a certain class of non-compact manifolds
(manifolds with a Lie structure at infinity \cite{ALN, LM, LN1,
MelroseCongress, MelroseScattering}).  This class of non-compact
manifolds is important for the analysis on locally symmetric spaces
\cite{Mueller1, Mueller2}.
  
Our main motivation to study gauge-invariant families comes from
analysis and spectral theory on non-compact manifolds.  For example,
the Fredholmness of an elliptic differential operator on a manifold
with a uniform structure at infinity is controlled by (possibly
several) gauge-invariant families of elliptic differential operators,
see \cite{LMN,LN0,LN1}.  This is a situation that arises when one
studies the Dirac operator on an $S^1$-manifold $M$, if we
desingularize the action of $S^1$ by replacing the original metric $g$
with $\phi^{-2}g$, where $\phi$ is the length of the infinitesimal
generator $X$ of the $S^1$-action. In this way, $X$ becomes of length
one in the new metric.  The main result of \cite{NistorDOP} states
that the kernel of the new Dirac operator on the open manifold $M
\smallsetminus M^{S^1}$ is naturally isomorphic to the kernel of the
original Dirac operator.  It turns out that the Fredholm property of
the resulting Dirac operator (obtained by the above procedure on the
non-compact manifold $M \smallsetminus M^{S^1}$) is controlled by the
invertibility of a family of operators invariant with respect to the
action of a bundle of non-abelian, solvable Lie groups.

In this paper, we study the gauge-equivariant index of families of 
elliptic operators invariant with respect to the free action a bundle of  
{\em simply-connected, solvable} Lie groups $\GR$. (By ``simply-connected'' we shall 
mean, as usual,  ``connected with trivial fundamental group.'')  This is the setting
derived from the example of the $S^1$-Dirac operator considered above
and is enough in many applications. The case of a bundle of compact 
Lie groups will be treated in a forthcoming paper with E. Troitsky
\cite{NistorTroi}. We are interested in the gauge-equivariant index 
of invariant families 
because the non-vanishing of the gauge-equivariant index is a first obstruction to
the invertibility of that family. One of the main results of this paper
is a formula for the Chern character of the gauge-equivariant index that  
is similar to the Atiyah-Singer index formula for families, but also
includes information on the {\em topology} of the bundle of Lie groups
that is considered, which is a qualitatively new feature of our setting.

We define the gauge-equivariant index of a family of invariant,  
elliptic operators using $K$-theory. Our definition works
for arbitrary family of Lie groups, proveded that they act freely.
(Both a {\em bundle of Lie groups} and a {\em family of Lie groups} are fiber
bundles with homeomorphic fibers. For a bundle of Lie groups
the fibers are required to be isomorphic as Lie groups. This
is not required, however, for families of Lie groups.)

In \cite{BismutCheeger}, Bismut and Cheeger have generalized the
Atiyah-Patodi-Singer index theorem \cite{APS} to families of Dirac
operators on manifolds with boundary (see also
\cite{Melrose-Piazza1}). Their results apply to operators whose
``indicial parts'' are invertible. These indicial parts are actually
families of Dirac operators invariant with respect to a one-parameter
group, so they fit into the framework of this paper (with $\GR = B
\times \RR$).  Thus, the results of this paper are relevant also to
the problem of studying families of operators on fibrations by
manifolds with boundary or fibrations by more general singular
spaces. (See also \cite{troitsky}.)
  
We now describe the contents of each section of this paper. In Section
\ref{S.Inv.Ops}, we discuss the action of a bundle of Lie groups $\GR
\to B$ on a fiber bundle $Y \to B$ and we introduce the algebras $\Psm
{\infty}Y$ of families of invariant pseudodifferential operators on
the fibers of $Y \to B$ whose convolution kernels are compactly
supported.  The classical definition of symbols applies to algebras
$\Psm{\infty}Y$ also.  In Section \ref{S.Abelian} and thereafter we
shall assume that $\GR$ acts freely on $Y$.  Section \ref{S.Abelian}
is devoted to the case when $\GR$ consists of simply-connected,
abelian Lie groups (which is the same thing as saying that $\GR \to B$
is a vector bundle, with the Lie group structure on each fiber given
simply by the addition of vectors). We prove then that the group of
gauge transformations of $\GR$ acts on $\Psm {\infty}Y$, and we use
this to enlarge the algebras $\Psm{\infty}Y$ by including families
whose convolution kernels are Schwartz functions. The resulting
algebras will be denoted by $\PsS {\infty}Y$.  Let $C_r^*(\GR)$ be the
closure of $\Psm{-\infty}Y$ in the topology defined by its action on
$L^2(Y)$. In Section \ref{S.Homotopy.Inv}, we consider arbitrary
families of Lie groups $\GR \to B$, but continue to assume that they
act freely on $Y$.  We define the gauge-equivariant index
\begin{equation*}
	\ind_\GR(A) \in K_0(C_r^*(\GR))
\end{equation*}
of a family of elliptic, invariant pseudodifferential operators $A \in
\Psm{\infty}Y$ (or $A \in \PsS{\infty}Y$, if $\GR$ is a vector
bundle). We then prove that the gauge-equivariant index of $A$ is the
obstruction to finding a regularizing family $R$ such that $A + R$,
acting between suitable Sobolev spaces, is invertible in each fiber,
provided that $\dim Y > \dim \GR$. This generalizes the usual property
of the Fredholm index of an elliptic pseudodifferential operator. In
Section \ref{Sec.Solvable}, we consider bundles of simply-connected,
solvable Lie groups $\GR \to B$, and we prove that
\begin{equation}\label{eq.K.solv}  
	K_*(C_r^*(\GR)) \simeq K^*(\lgg),  
\end{equation}  
provided that $\GR$ consists of simply-connected solvable Lie groups,
where $\lgg := \cup Lie\, \GR_b \to B$ is the vector bundle of Lie
algebras defined by $\GR$. Then we obtain a formula for the Chern
character of the gauge-equivariant index of an elliptic operator $A
\in \Psm {\infty}Y$.  (For simplicity, we called $A$ ``an operator,''
although it is really a family of operators. We shall do this
repeatedly.)  The proof of this result requires us to construct a
non-trivial deformation of the indicial family to a family to which
the Atiyah-Singer index theorem for families can be applied. We also
use here the notion of ``degree of an elliptic family.'' 
%which was
%first considered by Melrose.
  
Operators invariant with respect to $\RR^q$, a particular case of our
operators when the base is reduced to a point, appear in the
formulation of elliptic (or Fredholm) boundary conditions for
pseudodifferential operators on manifolds with corners. The
gauge-equivariant index can be used to study this problem, which is
relevant to the question of extending the Atiyah-Patodi-Singer
boundary conditions to manifolds with corners.  This is discussed in
the last section.
  
The algebras of invariant pseudodifferential operators that we study,
$\Psm {\infty}Y$ and $\PsS {\infty}Y$, are generalizations of
``parameter dependent'' algebras of pseudodifferential operators
considered by Agmon \cite{Agmon}, Grubb and Seeley \cite{GrubbSeeley},
Lesch and Pflaum \cite{LP}, Melrose \cite{MelroseScattering}, Shubin
\cite{SH}, and others. Our index theorem, Theorem \ref{Theorem.Chern},
thus gives a solution to the problem of determining the index of
elliptic, parameter dependent families of pseudodifferential
operators, parameterized by $\lambda \in \RR^q$. We note that the
concept of index of such a family requires a proper definition, and
that the Fredholm index, ``dimension of the kernel'' - ``dimension of
the cokernel,'' is not adequate for $q \ge 1$.  Actually, our
definition of the gauge-equivariant index is somewhat closer to the
definition of the $L^2$-index for covering spaces given in
\cite{Atiyah} and \cite{Singer} than to the definition of the Fredholm
index. However, an essential difference between their definition and
ours is that no trace is involved in our definition.
  
There exist several potential extensions of our results, although none
of them seems to be straightforward. These extensions will presumable
involve more general conditions on $\GR$ and its action on
$Y$. Probably the most general conditions under which one can
reasonably expect to obtain definite results are those when the fibers
of $\GR$ are connected and the action is proper. We need the fibers of
$\GR$ to be connected in order to be able to use general results on
the structure of connected Lie groups. The case when the fibers of
$\GR \to B$ are compact Lie groups will be discussed in
\cite{NistorTroi}, where the relevance of these construction to the
study of Ramond-Ramond charges and of the
Diximier-Douady-Donovan-Karoubi invariant.

We plan to apply the results of this paper to some problems in
$M$-theory \cite{SethiStern}, to the study of adiabatic limits of eta
invariants \cite{BF,Witten1}, and to analysis on locally-symmetric
spaces \cite{Mueller1,Mueller2}.

I would like to thank Fourier Institute and University of Grenoble,
where part of this project was completed, for their hospitality.
Also, I am grateful to Sergiu Moroianu and Evgenij Troitsky for useful
discussions. A preliminary version of this paper was circulated as a
Fourier Institute preprint in 1999.
  
All pseudodifferential operators considered in this paper are  
``classical,'' that is, one-step polyhomogeneous.

\section{Invariant families of (pseudo)differential operators\label{S.Inv.Ops}}

We now describe the settings in which we shall work, that of a family
of Lie groups $p: \GR \to B$ acting on a fiber bundle $\pi : Y \to B$.
Then we introduce the algebra $\Psm{\infty}Y$ consisting of families
of pseudodifferential operators along the fibers of $Y$ that are
invariant under the action of $\GR$ and whose convolution kernels are
compactly supported.  These algebras are particular cases of the
algebras introduced in \cite{NWX}. We can restrict in our discussion
to a connected component of $B$, so for simplicity, we shall assume
that $B$ is connected.  The gauge-equivariant index of the elliptic
operators $A \in \Psm{\infty}Y$ will form our main object of study.
  
Throughout this paper, $p:\GR \to B$ will denote a {\em family of Lie groups}  
on a smooth manifold $B$. By this we mean that $p: \GR \to B$ is  
a smooth fiber bundle, that each $\GR_b:= p^{-1}(b)$ is a Lie group, and that the  
multiplication and inverse depend differentiably on $b$. Hence the map  
\begin{equation}  
        \GR \times_B \GR := \{ (g',g) \in \GR \times \GR, d(g') = d(g)\}  
        \ni\, (g',g) \longrightarrow g'g^{-1}\, \in \GR  
\end{equation}  
is differentiable. This implies, by standard arguments, that the map  
sending a point $b \in \GR_b := p^{-1}(b)$ to $e_b$, the identity element of  
$\GR_b$, is a diffeomorphism onto a smooth submanifold of $\GR$. It  
also implies that the map $\GR \in g \to g^{-1} \in \GR$ is  
differentiable.  Note that we do not assume the groups $\GR_b$ to be
isomorphic, although this is true in most applications. 
Let $\lgg := \cup Lie\, \GR_b \to B$ be the bundle  
of Lie algebras of associated to $\GR$, that is, the bundle whose  
fiber above $b$, $\lgg_b$, is $Lie\, \GR_b$.  

Let us remark, although we shall not use this here, that $\GR$ is a 
particular case of a differentiable groupoid. The Lie algebroid associated
to this groupoid is the bundle $\lgg \to B$ defined above. (See \cite{NWX} or
\cite{LN1} for the necessary definitions.)
  
We also assume that $\GR$ acts smoothly on a fiber bundle $\pi: Y \to B$. 
This means that there are given actions $\GR_b \times Y_b \to Y_b := 
\pi^{-1}(b)$ of $\GR_b$ on $Y_b$, for each $b$, such that the induced map,  
$$  
        \GR \times_B Y :=\{ (g,y) \in \GR \times Y, d(g) = \pi(y)\}  
        \ni\, (g,y) \longrightarrow gy\, \in Y,  
$$  
is differentiable. The action of $\GR$ on $Y$ is called {\em free} 
[respectively, {\em proper}], if the action of $\GR_b$ on $Y_b$ is  
free [respectively, proper] for each $b$. (In this section we make no
assumptions on the action of $\GR$ on $Y$, other than that it is smooth,
but beginning with the next section we shall assume that this action 
is free.)
  
On $Y$, we consider smooth families $A = (A_b)$, $b \in B$, of
classical pseudodifferential operators acting on the fibers of $Y \to
B$. Unless mentioned otherwise, we assume that these operators act on
half densities along each fiber.  The family $(A_b)_{b \in B}$ will be
called {\em $\GR$-invariant} if each $A_b$ is invariant with respect
to the action of the group $\GR_b$. A {\em gauge-invariant family} is
a $\GR$-invariant family, for some family of Lie groups $\GR$.  The
algebra that we are interested in consists of $\GR$-invariant
operators satisfying also the following {\em support condition}. To
state this support condition, notice first that a family $A = (A_b)$
defines a continuous map $\CIc(Y) \to \CI(Y)$, and, as such, it has a
distribution (or Schwartz) kernel, which is a distribution $K_A$ on $Y
\times_B Y \subset Y \times Y$.  (We ignore the vector bundles in
which this distribution takes its values.)  Because the family
$A=(A_b)$ is invariant, the distribution $K_A$ is also invariant with
respect to the action of $\GR$. Consequently, the distribution $K_A$
is the pull back of a distribution $k_A$ defined on $(Y \times_B
Y)/\GR$. {\em We will require that $k_A$ have compact support.} We
shall sometimes call $k_A$ the {\em convolution kernel} of $A$.  This
condition on the support of $k_A$ ensures that each $A_b$ is a
properly supported pseudodifferential operator, and hence it maps
compactly supported functions (or sections of a vector bundle, if we
consider operators acting on sections of a smooth vector bundle) to
compactly supported functions (or sections). This support condition is
automatically satisfied if $Y/\GR$ is compact and each $A_b$ is a
differential operator. The space of smooth, invariant families $A$ of
order $m$ pseudodifferential operators acting on the fibers of $Y \to
B$ such that $k_A$ has compact support will be denoted by $\Psm{m}Y$.
Then
$$  
	\Psm {\infty}Y := \cup_{m \in \ZZ} \Psm{m}Y
$$   
is an algebra, by classical results
\cite{Hormander3,Taylor,TaylorB}. Note also that $\Psm{m}Y$ makes
sense also for $m$ not an integer.
  
We now discuss the principal symbols of an operator $A \in
\Psm{m}Y$. Let
$$  
        T_{vert}Y:=\ker (TY \to TB)  
$$   
be the bundle of {\em vertical} tangent vectors to $Y$, and let
$T_{vert}^*Y$ be its dual. We fix compatible metrics on $T_{vert}Y$
and $T_{vert}^*Y$, and define $S_{vert}^*Y$, the {\em cosphere bundle
of the vertical tangent bundle} to $Y$, to be the set of vectors of
length one of $T_{vert}^*Y$. Also, let
$$  
        \sigma_m : \Psi^m(Y_b) \to \CI(S_{vert}^* \cap T^*Y_b)  
$$   
be the usual principal symbol map, defined on the space of
pseudodifferential operators of order $m$ on $Y_b$. The definition of
$\sigma_m$ depends on the choice of a trivialization of the bundle of
homogeneous functions of order $m$ on $T^*_{vert}Y$, regarded as a
bundle over $S_{vert}^*Y$.  The principal symbols $\sigma_m(A_b)$ of
an element (or family) $A=(A_b) \in \Psm {m}Y$ then gives rise to a
smooth function on $\CI(S_{vert}^*Y)$, which is invariant with respect
to $\GR$, and hence descends to a smooth function on $S_{vert}^*Y$.
Let $\CI(S_{vert}^*Y)^\GR$ be the space of smooth functions on 
$S_{vert}^*Y$ that are invariant with respect to the canonical 
action of $\GR$. The resulting function,
\begin{equation}  
        \sigma_m(A) \in \CI(S_{vert}^*Y)^\GR,  
\end{equation}  
will be referred to as the {\em principal symbol} of an element (or
operator) in $\Psm {m}Y$.
  
In the particular case $Y = \GR$, $\Psm {\infty}\GR$ identifies with
families of convolution operators on the fibers $\GR_b$ with kernels
contained in a compact subset of $\GR$, smooth outside the identity,
and only with conormal singularities at the identity. In particular,
$\Psm {-\infty}{\GR}= \CIc(\GR)$, with the fiberwise convolution
product.
  
Suppose now that the quotient $Y/\GR$ is compact, which implies that
$(S^*_{vert}Y)/\GR$ is also compact. As is customary, an operator $A
\in \Psm {m}Y$ is called {\em elliptic} if, and only if, its principal
symbol is everywhere invertible. The same definition applies to $A =
[A_{ij}] \in M_N(\Psm {m}Y)$. Namely, the operator $A$, regarded as
acting on sections of the trivial vector bundle $\CC^N$, is elliptic
if, and only if, its principal symbol
$$   
        \sigma_m(A) := [\sigma_m(A_{ij})] \in
        M_n(\CI(S_{vert}^*Y)^\GR)
$$    
is invertible.  
  
Assume that there is given a $\GR$-invariant metric on $T_{vert}Y$,
the bundle of {vertical} tangent vectors, and a $\GR$-equivariant
bundle $W$ of modules over the Clifford algebras of $T_{vert}Y$. Then
a typical example of a family $D=(D_b) \in \Psm {\infty}Y$ is that of
the family of Dirac operators $D_b$ acting on the fibers $Y_b$ of $Y
\to B$.  (Each $D_b$ acts on sections of $W\vert_{Y_b}$, the
restriction of the given Clifford module $W$ to that fiber.)

The considerations of this section extend immediately to operators
acting between sections of two $\GR$-equivariant vector bundles. If
$E_0$ and $E_1$ are $\GR$-equivariant vector bundles, we denote by
$\Psm{m}{Y;E_0,E_1}$ the space of $\GR$-invariant families of order
$m$ pseudodifferential operators acting on sections of $E_0$, with
values sections of $E_1$, whose convolution kernels have compact
support.

\section{The case when $\GR$ is a vector bundle \label{S.Abelian}}

{\em Beginning with this section and throughout the rest of this
paper, we shall assume that the action of $\GR \to B$ on $Y \to B$ is
free} (that is, that the action of $\GR_b$ on $Y_b$ is free, for any
$b \in B$). {\em In this section we shall also assume that the fibers
of $\GR \to B$ are simply-connected, abelian Lie groups.}  This is
assumption is equivalent to the assumption that $\GR \to B$ is a
vector bundle with the induced bundle of Lie algebras structure.

For the discussion of the gauge-equivariant index of a an operator $A
\in \Psm{m}Y$ in Section \ref{S.Homotopy.Inv}, we need to introduce
the algebra $\PsS{\infty}Y$, which is a variant of the algebra
$\Psm{\infty}Y$ which it contains. The reason for considering the
slightly larger algebra $\PsS{\infty}Y$ is that it will have the
following property, if $A = (A_b) \in \PsS{\infty}Y$ and each $A_b$ is
invertible on $L^2(Y_b)$ (in the sense of unbounded operators), then
$(A_b^{-1}) \in \PsS {\infty}Y$. (This property of the algebra
$\PsS{\infty}Y$ is usually referred to in the literature as ``spectral
invariance.'')

Take $ Y = B \times Z \times \RR^q, $ with $Z$ a compact manifold and
$\GR = B \times \RR^q, $ $\pi : Y \to B$ and $p : \GR \to B$ being the
projections onto the first components of each product.  The action of
$\GR$ on $Y$ is given by translation on the last component of
$Y$. Then the $\GR$-invariance condition becomes simply $\RR^q$
invariance with respect to the resulting $\RR^q$ action. If $Y$ and
$\GR$ are as described here, then we call $Y$ {\em a flat
$\GR$-space}.  The construction of the algebra $\PsS{\infty}Y$ is
local, so we may assume that $Y$ is a flat $\GR$-space.
     
The residual ideal of the algebra $\Psm {\infty}Y$ is $\Psm
{-\infty}Y$ and consists of operators that are regularizing along each
fiber. More precisely, it consists of those families of smoothing
operators on $Y = B \times Z \times\RR^q$ that are
translation-invariant under the action of $\RR^q$ and have {\em
compactly supported} convolution kernels. Thus
$$  
        \Psm {-\infty}Y \cong \CIc(B \times \RR^q; \Psi^{-\infty}(Z))
        \subset \mathcal{S}(B \times \RR^q; \Psi^{-\infty}(Z)) \simeq
        \mathcal{S}(B \times Z\times Z\times\RR^q).
$$  
(Here ${\mathcal S}$ is the generic notation for the space of Schwartz
functions.) The second isomorphism above is obtained from the
isomorphism
$$  
        \Psi^{-\infty}(Z) \simeq \CI(Z \times Z),  
$$  
defined by the choice of a nowhere vanishing density on $Z$. If we
also endow $\Psi^{-\infty}(Z)$ with the locally convex topology
induced by this isomorphism, then it becomes a nuclear, locally convex
space.  We now enlarge the algebra $\Psm {\infty}Y$ to include all
$\GR$-invariant, regularizing operators whose kernels are in
$\mathcal{S}(B \times Z\times Z\times\RR^q)$:
\begin{equation}  
        \PsS{m}Y = \Psm{m}Y + \mathcal{S}(B \times Z\times  
        Z\times\RR^q).  
\end{equation}   
Explicitly, the action of $T \in \mathcal{S}(B \times Z\times  
Z\times\RR^q)$ on a smooth function $f \in \CIc(B \times Z \times  
\RR^q)$ is given by  
$$  
        Tf(b,y_1,t) = \int_{Z \times \RR^q} T(b, y_1,  
        y_0,t-s)f(b,y_0,s) dy_0 ds.  
$$  
For $B$ reduced to a point, the algebras $\PsS {\infty}Y$ were
considered before (see \cite{MelroseScattering} and the references
therein) as the range space of the indicial map for ``cusp''
pseudodifferential operators on manifolds with corners
  
Still assuming that we are in the case of a flat $\GR$-space, we
notice that the Fourier transformation ${\mathcal F}_2$ in the
translation-invariant directions gives a dual identification
\begin{equation}  
        \PsS {-\infty}Y \ni T \to \hat{T}:=  
        {\mathcal F}_2 T {\mathcal F}_2^{-1}  
        \in \mathcal{S}(B \times Z\times Z \times \RR^q) =  
        \mathcal{S}(B \times \RR^q; \Psi^{-\infty}(Z))  
\label{Smoothing.Ideal}  
\end{equation}  
with the space of Schwartz functions with values in the smoothing  
ideal of $Z.$ The point of this identification is that the  
convolution product is transformed into the pointwise product. The  
Schwartz topology from \eqref{Smoothing.Ideal} then gives $\PsS  
{-\infty}Y$ the structure of a nuclear locally convex topological  
algebra. Following the same recipe, the Fourier transform also gives  
rise to an {\em indicial family}  
\begin{equation}\label{Fourier.Transform}  
        \Phi:\PsS {\infty}Y\ni T \to \hat{T}:={\mathcal F}_2 T
        {\mathcal F}_2^{-1} \in \CI(B \times \RR^q;\Ps{\infty}{Z}).
\end{equation}  
We denote   
$ 
        \hat{T}(\tau)=\Phi(T)(\tau).  
$
  
The map $\Phi$ is not an isomorphism since $\hat{A}(\tau)$ has joint
symbolic properties in the variables of $\RR^q$ and the fiber
variables of $T^*Z.$ Actually, the principal symbols of the operators
$\hat{A}(\tau)$, for $\tau$ in a fixed fiber of $T_{vert}^*Y \to B$,
is independent of $\tau$.

\begin{lemma}\label{Lemma.dilations}\ Assume $Y = B \times Z \times \RR^q$   
is a flat $\GR$-space. Then the action of the group $GL_q(\RR)$ on the
last factor of $Y = B \times Z \times \RR^q$ extends to an action by
automorphisms of $\CI(B,GL_q(\RR))$ on $\Psm {\infty} Y$ and on $\PsS
{\infty} Y.$
\end{lemma}

\begin{proof}\   
The vector representation of $GL_q(\RR)$ on the second component of $Z
\times \RR^q$ defines an action of $GL_q(\RR)$ on $\Psi^\infty(Z
\times \RR^q)$ that preserves the class of properly supported
operators and the products of such operators.  It also normalizes the
group $\RR^q$ of translations, and hence it maps $\RR^q$-invariant
operators to $\RR^q$-invariant operators. This property extends right
away to the action of $\CI(B,GL_q(\RR))$ on families of operators on
$B \times Z \times \RR^q$, and hence $\CI(B,GL_q(\RR))$ maps $\Psm
{m}Y$ isomorphically to itself. From the isomorphism
\eqref{Smoothing.Ideal}, we see that $\CI(B,GL_q(\RR))$ also maps
$\PsS {-\infty}Y$ to itself. This gives an action by automorphisms of
$\CI(B,GL_q(\RR))$ on $\PsS {\infty} Y$, which is the sum of $\PsS
{-\infty}Y$ and $\Psm {\infty}Y$.
\end{proof}

By choosing a lift of $Y /\GR \to Y$, which is possible because the
fibers are contractable ($\GR$ is a vector bundle), we obtain that,
locally, the bundle $Y$ is isomorphic to a flat $\GR$ space.  Then the
above lemma allows us to extend the previous definitions, including
those of the algebras $\PsS {\infty}Y$ and of the indicial family from
the flat case to the case $\GR$ a vector bundle. Recall the bundle
$\mfk g : = \cup Lie\, \GR_b \to B$.  The indicial family $\hat{A}$ of
an operator $A \in \PsS {\infty}Y$, will then be a family of
pseudodifferential operators acting on the fibers of $Y \times_B
\lgg^* \to \lgg^*$ (here $\lgg^*$ is the dual of the vector bundle
$\lgg$):
\begin{equation}\label{eq.bundle.hatA}   
        \hat{A}(\tau) \in \Psi^*(Y_b/\GR_b), \quad \text{if } \tau \in  
        \lgg_b^* := (Lie\, \GR_b)^*.  
\end{equation}   
The action of $GL_q(\RR)$ in the above lemma will have to be replaced with   
the group of gauge-transformations of $\lgg$.  
  
We can look at a general vector bundle $\GR$ from an other, related  
point of view. Let $Z' := Y/\GR$, which is again a fiber bundle $Z' \to B$. 
The algebra $\Psi^{\infty}(Z' \vert B)$ of smooth families of pseudodifferential   
operators along the fibers of $Z' \to B$ can be regarded as the algebra of  
sections of a vector bundle on $Z'$ (with infinite dimensional fibers).  
We can lift the vector bundle associated $\Psi^\infty (Z' \vert B)$ to $\lgg^*$ 
and then $\hat{A}$ is a section of this lifted vector bundle over $\lgg^*$. 
We shall write this as follows:   
\begin{equation}\label{eq.section.hatA}   
	\hat{A} \in \Gamma(\lgg^*,\Psi^{\infty}(Z' \vert B)).  
\end{equation}

\section{The gauge-equivariant index: Definition\label{S.Homotopy.Inv}}

In this section, $\GR \to B$ will an arbitrary family of Lie groups
(we shall make no assumptions on the fibers $\GR_b$), but we shall
continue to assume that the action of $\GR$ on $Y$ is free

We now define study $\GR$-invariants of elliptic operators in
$M_N(\Psm {\infty}Y)$, the main invariant being the gauge-equivariant
index of such an invariant, elliptic family. If $\dim Y > \dim \GR$,
we then show that the gauge-equivariant index gives the obstruction
for family $A=(A_b) \in M_N(\Psm{m}Y)$ to have a perturbation by a
family $R=(R_b) \in M_N(\Psm{-\infty}Y)$, such that $A+R = (A_b +
R_b)$ be invertible, for all $b \in B$, between suitable Sobolev
spaces, see Theorem \ref{Theorem.obst}. For families of
simply-connected, abelian Lie groups $\GR$, we give an interpretation
of the gauge-equivariant index of an elliptic operator in terms of its
indicial family. This requires the notion of ``degree,'' due to
Melrose. We shall use this in the next section to prove an
Atiyah-Singer index type formula for the Chern character of the index
of a family of $\GR$-invariant, elliptic operators. If $\GR$ is a
vector bundle, then we can consider the algebra $\PsS{\infty}Y$
instead of $\Psm {\infty}Y$.
  
We now proceed to define the gauge-equivariant index of an elliptic
family $A \in \Psm {m}Y$.  This will be done using the $K$-theory of
Banach algebras. Let $C^*_r(Y,\GR)$ be the closure of $\Psm
{-\infty}{Y}$ with respect to the norm
$$  
        \| A \| = \sup_{b \in B} \|A_b\|  
$$  
each operator $A_b$ acting on the Hilbert space of square integrable
densities on the fiber $Y_b$. If $Y = \GR$, then we also write
$C^*_r(\GR,\GR)=C^*_r(\GR)$. (We remark, for those familiar with the
concepts, that $C^*_r(\GR)$ is the reduced $C^*$-algebra associated to
the groupoid $\GR$. Also, the full and the reduced $C^*$-algebras,
$C^*(\GR)$ and $C_r^*(\GR)$ coincide, if the fibers $\GR_b$ are
solvable groups.) For each locally compact space $X$, we denote by
$C_0(X)$ the space of continuous functions vanishing at infinity on
$X$.  Then, if $\GR$ is a vector bundle, we have $C^*_r(\GR) \simeq
C_0(\GR^*)$.
  
We shall use below $\widehat{\otimes}$, the minimal tensor product of  
$C^*$-algebras. This minimal tensor product is defined to be  
(isomorphic to) the completion of the image of $\pi_1 \otimes \pi_2$,  
the tensor product of two injective representations $\pi_1$ and  
$\pi_2$.  For the cases we are interested in, the minimal and the  
maximal tensor product coincide \cite{Sakai}.

\begin{lemma}\label{lemma.stable}\   
Assume as above that $\GR$ acts freely on $Y$ and that
$\dim Y > \dim \GR$. Also, let ${\mathcal K}={\mathcal  
K}(Y_b/\GR_b)$ denote the algebra of compact operators on one of the  
fibers $Y_b/\GR_b$, for some fixed but arbitrary $b \in B$. Then  
$$  
        C_r^*(Y,\GR) \simeq C_r^*(\GR) \widehat{\otimes} {\mathcal K}.  
$$  
Consequently, $K_i(C_r^*(Y,\GR)) \simeq K_i(C_r^*(\GR))$.  
\end{lemma}

\begin{proof}\ Let $\mfk A$ be the space of sections of  
the bundle of algebras ${\mathcal K}(Y_b/\GR_b)$.  If $Y$ is a flat  
$\GR$-space, then the isomorphism $C_r^*(Y,\GR) \simeq C_r^*(\GR)  
\widehat{\otimes} {\mathcal K}$ follows, for example, from the results  
of \cite{LMN}. In general, this local isomorphism gives 
$C_r^*(Y,\GR) \simeq C_r^*(\GR) \widehat{\otimes}_{C_0(B)} {\mathfrak 
A}$. 
  
Our assumptions imply that ${\mathcal K}$ is infinite dimensional, and  
hence its group of automorphisms is contractable, see  
\cite{Diximier}. Consequently, there is no obstruction to trivialize  
the bundle of algebras ${\mathcal K}(Y_b/\GR_b)$, and hence  
${\mathfrak A} \simeq C_0(B) \widehat{\otimes}{\mathcal K}$. We obtain  
\begin{equation*}   
        C_r^*(Y,\GR) \simeq C_r^*(\GR) \widehat{\otimes}_{C_0(B)} 
        {\mathfrak A} \simeq C_r^*(\GR) \widehat{\otimes}_{C_0(B)} 
        \big[C_0(B) \widehat{\otimes}{\mathcal K}\big] \simeq C_r^*(\GR) 
        \widehat{\otimes} {\mathcal K}. 
\end{equation*}    
   
The last part of lemma follows from the above isomorphisms and from
the natural isomorphism $K_i(A \widehat{\otimes} {\mathcal K}) \simeq
K_i(A)$, valid for any $C^*$-algebra $A$.
\end{proof}

We proceed now to define the gauge-equivariant index of an {\em elliptic, $\GR$-invariant  
family} of operators  
$$   
        A = (A_b) \in M_N(\Psm{m}Y) = \Psm{m}{Y;\CC^N}  
$$   
when $Y/\GR$ is compact. More general operators can be handled similarly,
and the necessary changes will discussed below. We observe first that there 
exists an exact sequence  
\begin{equation}\label{eq.exact.seq}  
  	0 \to C^*_r(Y,\GR) \to {\mathcal E} \to \CI(S^*_{vert}Y) \to  
  	0,\quad {\mathcal E}:= \Psm {0}Y + C^*_r(Y,\GR),  
\end{equation}  
obtained using the results of \cite{LMN}. The operator $A$  
(or, rather, the family of operators $A = (A_b)$) has an invertible  
principal symbol, and hence the family $T = (T_b)$,  
$$   
        T_b:=(1 + A_b^*A_b)^{-1/2}A_b,  
$$  
consists of elliptic, $\GR$-invariant operators. Moreover, the
operator $T$ is an element of ${\mathcal E} := \Psm {0}Y +
C^*_r(Y,\GR)$, its principal symbol is invertible, and hence it
defines a class $[T] \in K_1(\CI(S^*_{vert}Y)) \simeq
K^1(S^*_{vert}Y)$. Let
\begin{equation*}\label{eq.boundary.map}  
        \pa : K_1^{alg}(S^*_{vert}Y) \to K_0^{alg}(C^*_r(Y,\GR))
        \simeq K_0 (C^*_r(Y,\GR))
\end{equation*}  
be the boundary map in the $K$-theory exact sequence  
\begin{multline*}  
  	K_1^{alg}(C^*_r(Y,\GR)) \to K_1^{alg}({\mathcal E}) \to
  	K_1^{alg}(S^*_{vert}Y)
  	\stackrel{\pa}{\rightarrow}K_0^{alg}(C^*_r(Y,\GR)) \\ \to
  	K_0^{alg}({\mathcal E}) \to K_0^{alg}(S^*_{vert}Y)
\end{multline*}  
associated to the exact sequence \eqref{eq.exact.seq}. Because $K_0 
(C^*_r(Y,\GR)) \simeq K_0(C^*_r(\GR)),$ by Lemma \ref{lemma.stable}, 
we get a group morphism 
\begin{equation}  
        \ind_\GR : K_1^{\alg}(\CI(S^*_{vert}Y))  \to  K_0(C^*_r(\GR)),  
\end{equation}  
which we shall call {\em the gauge-equivariant index morphism}. The
image of $A$ under the composition of the above maps is
$\ind_\GR([T])$. We shall also write $\ind_\GR(A) := \ind_\GR([T])$,
and call $\ind_\GR(A)$ {\em the gauge-equivariant index of} $A$. A
more direct but longer definition is contained in the proof of Theorem
\ref{Theorem.obst} below (see Equation \eqref{eq.def.ind}).  The
gauge-equivariant index morphism descends in this case to a group
morphism $K_1^{\alg}(\CI(S^*_{vert}Y)) \to K_0(C^*_r(\GR))$ denoted in
the same way.  If $\GR$ is a vector bundle, we can replace ${\mathcal
E}$ with $\PsS {0}Y$ and $C^*_r(\GR)$ with $\PsS {-1}{Y}$.
  
We denote by $I_N$ the unit of the matrix algebra $M_N(\mathcal E)$.  
If $A \in \Psm{0}{Y;E}$, we can find $N$ large such that $\Psm{0}{Y;E}  
\subset M_N(\Psm{0}Y)$ with $1$ mapping to the projector $e$ under  
this isomorphism. Then $A + (I_N-e)$ is invertible in this matrix  
algebra, and we define then $\ind_\GR(A) = \ind_\GR(A + I_N - e)$.  

Let us briefly see what needs to be changed when $Y /\GR$ is not
compact.  If $Y/\GR$ is not compact, the above definition of the
gauge-equivariant index applies only to elliptic operators in $I_N +
M_N(\Psm{0}{Y;E})$. However, all results below extend to operators in
$I_N + M_N(\Psm{0}{Y;E})$, after some obvious changes.  For
differential operators acting between sections of {\em different}
bundles, we can define the gauge-equivariant index using the adiabatic
groupoid of $\GR$ as in \cite{LMN}. For arbitrary elliptic
pseudodifferential operators $A \in \Psm{m}{Y;E_0,E_1}$, acting
between sections of possibly different vector bundles, we can define
the gauge-equivariant index using Kasparov's bivariant
$K$-Theory. Since, when $Y/\GR$ is non-compact, no element in
$\Psm{m}{Y;E_0,E_1}$ is invertible modulo regularizing operators, we
must allow more general operators in this case. For example, we can
take
\begin{equation*}
	A \in \Gamma(Y,\Hom(E_0,E_1))^{\GR} + \Psm{0}{Y;E_0,E_1},
\end{equation*}
provided that it is bounded. If $A \in \Psm{m}{Y;E_0,E_1}$, $m >0$, we
replace $A$ by $(1 + AA^*)^{-1/2} A$, which will be an operator in the
closure of $\Psm{0}{Y;E_0,E_1}$, by the results of \cite{LN1}. If $A$
is elliptic, then the resulting family provides directly from the
definition an element
\begin{equation}
	\ind_\GR(A) \in KK(\CC, C_r^*(\GR)).
\end{equation}
To obtain the gauge-equivariant index as an element of
$K_0(C_r^*(\GR))$, we use the isomorphism $KK(\CC, C_r^*(\GR)) \cong
K_0(C_r^*(\GR))$. We can also work with $\PsS{-\infty}{\GR}$ instead
of $C_r^*(\GR)$, by using the smooth versions of $KK$-theory
introduced in \cite{NistorKthry,NistorANN}.
  
The main property of the gauge-equivariant index of an operator $A$ is
that it gives the obstruction to the existence of invertible
perturbations of $A$ by lower order operators. We denote by $H^s(Y_b)$
the $s$th Sobolev space of $1/2$-densities on $Y_b$, which is uniquely
defined because of the bounded geometry of $Y_b$, for $Y/\GR$
compact. More precisely, $H^s(Y_b)$ is, by definition, the domain of
$(1 + D^*D)^{s/2m}$, if $D \in \Psm{m}Y$ is elliptic and $s \ge
0$. For $s \le 0$, $H^s(Y_b)$ is, by definition, the dual of
$H^{-s}(Y_b)$. If $E \to Y$ is a $\GR$-equivariant vector bundle, then
we shall denote by $H^s(Y_b,E)$ the Sobolev space of sections of $E$,
defined analogously.

\begin{theorem}\label{Theorem.obst}\   
Let $\GR \to B$ be a bundle of Lie groups acting on the fiber bundle
$Y \to B$, as above, and assume that $Y/\GR$ is compact, of positive
dimension.  Let $A \in \Psm {m}{Y,E_0,E_1}$ be an elliptic operator.
Then we can find $R \in \Psm {m-1}{Y,E_0,E_1}$ such that
$$   
        A_b + R_b : H^{s}(Y_b,E_0) \to H^{s-m}(Y_b,E_1)   
$$    
is invertible for all $b \in B$ if, and only if, $\ind_\GR(A)=0$.  
Moreover, if $\ind_\GR(A)=0$, then we can choose $R \in \Psm{-\infty}Y$.  
The same result holds if $Y/\GR$ is non-compact and  
$A \in \Gamma(Y,\Hom(E_0,E_1))^{\GR} + \Psm {0} {Y;E_0,E_1}$ is elliptic  
and bounded.   
\end{theorem}

\begin{proof}\  
We shall assume that $E_0=E_1$ are trivial vector bundles of rank $n$,
for simplicity of notation. The proof in the general case is the same.
It is clear from definition that if we can find $R$ with the desired
properties, then $\ind_\GR(A) = 0 \in K^0(\lgg^*)$.  Suppose now that
$A \in \Psm {m}{Y,\CC^N}$ is elliptic and has vanishing
gauge-equivariant index.  Using the notation $T = (1 + AA^*)^{-1/2}A$,
we see that $A_b$ is invertible between the indicated Sobolev spaces
if, and only if, $T_b$ is invertible as a bounded operator on
$L^2(Y_b)^N$. We can hence assume that $m = 0$ and $T = A$.
  
Because $C^*_r(Y,\GR)$ satisfies $C^*_r(Y,\GR)\simeq C^*_r(\GR)
\widehat{\otimes} {\mathcal K}$, by Lemma \ref{lemma.stable}, we can
use some general techniques to prove that the vanishing of
$\ind_\GR(A)$ implies that $A$ has a perturbation by $\GR$-invariant,
regularizing operators in $\Psm{-\infty}{Y,\CC^N}$ that is invertible
on each fiber. We fix an isomorphism $M_N(C^*_r(Y,\GR)) \cong
C^*_r(\GR) \otimes {\mathcal K}$. We now review this general technique
using a generalization of an argument from \cite{NistorKthry}. Let
$\mathcal E$ be the algebra introduced in Equation
\eqref{eq.exact.seq}.  We denote by $I_N$ the unit of the matrix
algebra $M_N(\mathcal E)$.  Also, denote by $\overline{\mathcal E}$
the closure of $\mathcal E$ in norm.
  
Choose a sequence of projections $p_n \in {\mathcal K}$, $\dim p_n =
n$, such that $p_n \to 1$ in the strong topology. Because $A_b$ is
invertible modulo $C^*_r(\GR_b) \otimes {\mathcal K}$, we can find a
large $n$ and $R \in M_N(\Psm{-\infty}Y)$ such that
$$  
	A'_b := A_b \oplus R_b : L^2(Y_b)^N \oplus L^2(Y_b)^N \to L^2(Y_b)^N  
$$   
is surjective and $(1 \otimes p_n) R_b (1 \otimes p_n) = R_b$, for all
$b \in B$. Then $\{0\}$ is not in the spectrum of ${A'}{A'}^*$, and we
can consider $V := (A'{A'}^*)^{-1/2} A' \in M_N(\overline{\mathcal
E})$, which by construction satisfies $VV^* = I_N \in
M_N(\overline{\mathcal E})$. Consequently, $V^*V$ is a projection in
$M_{2N}(\overline{\mathcal E})$. Because $A$ is invertible modulo
$M_N(C^*_r(\GR) \otimes {\mathcal K})$,
$$  
	V^*V - VV^* \in  M_{2N}(C^*_r(\GR) \otimes {\mathcal K}).  
$$  
Let $e = I_N \oplus (1 \otimes p_n) -V^*V$, which is also a
projection, by construction. Moreover,
\begin{equation}\label{eq.def.K}  
	e - (1 \otimes p_n) \in M_{2N}(C^*_r(\GR) \otimes \mathcal K),  
\end{equation}  
and hence both $e$ and $1 \otimes p_n$ are projection in   
$M_{2N}(C^*_r(\GR)^+ \otimes \mathcal K)$ (for any algebra $B$,   
we denote by $B^+$ the algebra with adjoint unit). Equation \eqref{eq.def.K}  
gives that, by definition, $[e] - [1 \otimes p_n]$ defines a $K$-theory class  
in $K_0(C^*_r(\GR))$. From definition, we get then  
\begin{equation}\label{eq.def.ind}  
	\ind_\GR(A) = [e] - [1 \otimes p_n].  
\end{equation}  
  
Now, if $\ind_\GR(A)=0$, then we can find $k$ such that $e \oplus
1\otimes p_k$ is (Murray-von Neumann) equivalent to $1\otimes p_{n +
k}$. By replacing our original choice for $n$ with $n + k$, we may
assume that $e$ and $1 \otimes p_n$ are equivalent, and hence that we
can find $u \in \CC I_{2N} + M_{2N}(C^*_r(\GR) \otimes {\mathcal K})$
with the following properties:\ there exists a large $l$ and $x \in
M_{2N}(C^*_r(\GR))$ such that, if we denote $e_0 = 1 \otimes p_l
\oplus 1 \otimes p_n$, then $u = I_{2N} + e_0 x e_0$ and $e = u(1
\otimes p_n) u^{-1}$. Then $Vu$ is in $M_n (\overline{\mathcal E})$
(more precisely $I_N V I_N = V$) and is invertible.  Consequently
$B:=(A'{A'}^*)^{1/2} V$ is also invertible. But $B$ is a perturbation
of $A'$, and hence also of $A$, by an element in $M_{2N}(C^*_r(\GR)
\otimes {\mathcal K})$. Since $\Psm{-\infty}Y$ is dense in
$C^*_r(\GR)$, this gives the result.
\end{proof}

\section{The gauge-equivariant index: A formula\label{Sec.Solvable}}

We now treat in more detail the case of bundles of solvable Lie  
groups, when more precise results can be obtained. Other classes of  
groups will lead to completely different problems and results, so we  
leave their study for the future. The class of simply-connected  
solvable fibers is rich enough and has many specific features, so we  
content ourselves from now on with this case only.  \\[2mm]  
{\bf Assumption.}\ {\em From now on and throughout this paper, we
shall assume that the family $\GR$ consists of simply-connected
solvable Lie groups.} By ``simply-connected'' we mean, as usual,
``connected with trivial fundamental group.''
  
\vspace*{2mm} We shall denote by $\lgg \to B$ the bundle of Lie  
algebras of the groups $\GR_b$, $\lgg_b \simeq Lie\,\GR_b$ and by  
$$  
        \exp: \lgg \to \GR.  
$$  
the exponential map.   
    
In order to study the algebra $C^*_r(\GR)$ and its $K$-groups, we
shall deform it to a commutative algebra. This deformation is obtained
as follows.  Let $\GR_{adb} = \{0\} \times \lgg \cup (0,1] \times
\GR$, $B_1 = [0,1] \times B$, and $p_{adb}:\GR_{adb} \to B_1$ be the
natural projection. On $\GR_{adb}$ we put the smooth structure induced
by
$$  
        \phi: B_1 \times \lgg \to \GR_{adb}  
$$  
$\phi(0,X) = (0,X)$ and $\phi(t,X)=(t,\exp(tY))$, which is a bijection
for all $(t,X) \in [0,1] \times \lgg$ in a small neighborhood of $\{0
\} \times \lgg \cup B_1 \times \{0\}$.  Then we endow $\GR_{adb}$ with
the Lie bundle structure induced by the pointwise product.  Evaluation
at $t \in [0,1]$ induces algebra morphisms
$$  
        e_t : \Psm {m}{\GR_{adb}} \to \Psm {m}{\GR}, \quad t > 0,  
$$  
and   
$$  
        e_0 : \Psm {m}{\GR_{adb}} \to \Psm {m}{\lgg}, \quad t = 0.  
$$  
Passing to completions, we obtain morphisms $e_t$ from
$C_r^*(\GR_{adb})$ to $C_r^*(\GR)$, for $t > 0$, and to
$C_r^*(d^{-1}(0) \times B) \simeq C_0(\lgg^*)$, for $t = 0$.

\begin{lemma}\label{Lemma.Hom.Inv}\   
The morphisms $e_t : C_r^*(\GR_{adb}) \to C_r^*(\GR)$, for $t > 0$, and  
$e_0 : C_r^*(\GR_{adb}) \to C_0(\lgg^*)$, for $t = 0$, induce  
isomorphisms in $K$-theory.  
\end{lemma}

\begin{proof}\ Assume first that there exists a Lie group bundle morphism   
$\GR \to B \times \RR$. (In other words, there exists a smooth map
$\GR \to \RR$ that is a morphism on each fiber.) Let $\GR'$ denote the
kernel of this morphism and let $\GR_{adb}'$ be obtained from $\GR'$
by the same deformation construction by which $\GR_{adb}$ was obtained
from $\GR$. Then we obtain a smooth map $\GR_{adb} \to \RR$ that is a
group morphism on each fiber, and hence
$$   
        C_r^*(\GR_{adb}) \simeq C_r^*(\GR_{adb}') \rtimes \RR, \;
        C_r^*(\GR) \simeq C_r^*(\GR') \rtimes \RR, \; \text{ and }
        C_0(\lgg^*) \simeq C_0({\lgg'}^*) \rtimes \RR.
$$   
Moreover, all above isomorphisms are natural, and hence compatible with   
the morphisms $e_t$. Assuming now that the result was proved for all   
Lie group bundles of smaller dimension, we obtain the desired result   
for $\GR$ using Connes' Thom isomorphism in $K$-theory \cite{ConnesThom},    
which in this particular case gives:    
$$   
        K_i(C_r^*(\GR_{adb})) \simeq K_{i+1} (C_r^*(\GR_{adb}')), \;
        K_i(C_r^*(\GR)) \simeq K_{i+1}(C_r^*(\GR')), \; \text{ and }
$$  
$$   
        K_i(C_0(\lgg^*)) \simeq K_{i+1} (C_0({\lgg'}^*)).   
$$   
   
This will allow us to complete the result in the following way. Let
$U_k$ be the open subset of $B$ consisting of those $b \in B$ such
that $[\GR_b,\GR_b]$ has dimension $\ge k$. From the Five Lemma and
the six term exact sequences in $K$-theory associated to the ideal
$C_r^*(\GR_{adb}\vert_{U_{k} \smallsetminus U_{k+1}})$ of
$C_r^*(\GR_{adb}\vert_{U_{k-1} \smallsetminus U_{k+1}})$, for each
$k$, we see that it is enough to prove our result for
$\GR_{adb}\vert_{U_k \smallsetminus U_{k+1}}$ for all $k$. Thus, by
replacing $B$ with $U_k \smallsetminus U_{k+1}$, we may assume that
the rank of the abelianization of $\GR_b$ is independent of $b$.
Consequently, the abelianizations of $\GR_b$ form a vector bundle
$$   
        {\mathcal A} := \cup \GR_b / [\GR_b,\GR_b]   
$$   
on $B$.   
   
A similar argument, using the Meyer-Vietoris exact sequence in
$K$-theory and the compatibility of $K$-theory with inductive limits
\cite{Blackadar}, shows that we may also assume the vector bundle
${\mathcal A}$ of abelianizations to be trivial. Then the argument at
the beginning of the proof applies, and the result is proved.
\end{proof}

{}From the above lemma we immediately obtain the following corollary:

\begin{corollary}\label{Cor.Hom.Inv}\   
Let $\GR$ be a bundle of simply connected, solvable Lie groups. Then  
$$  
        K_i(C_r^*(\GR)) \simeq K_i(C_0(\lgg^*)) \simeq K^i(\lgg^*).  
$$  
\end{corollary}

We now give an interpretation of $\ind_\GR(A)$, for $\GR$ a vector bundle, using  
the properties of the indicial family $\hat{A}(\tau)$ of $A$. We assume  
that $Y/\GR$ is compact.  
   
We shall also use the following construction.  Let $X$ be a compact  
manifold with boundary.  Let $T(x)$ be a family of elliptic  
pseudodifferential operators acting between sections of two smooth  
vector bundles, $E_0$ and $E_1$, on the fibers of a fiber bundle $M  
\to X$ whose fibers are compact manifolds without corners. Then we can  
realize the index of $T$ as an element in the relative group  
$K^0(X,\pa X)$. This can be done directly using Kasparov's theory or  
by the ``Atiyah-Singer trick'' as follows. We proceed as in  
\cite{AS4}, Proposition (2.2), to define a smooth family of maps  
$R(x): \CC^N \to \CI(Y)$, such that the induced map  
$$  
	V := T \oplus R : \CI(X)^N \oplus \CI(X,L^2(Y;E_0)) \to   
	\CI(X,L^2(Y;E_1))  
$$  
is onto for each $x$. Since $T(x)$ is invertible for $x \in \pa X$, we  
can choose $R(x)=0$ for $x \in \pa X$. Then $\ker (V)$ is a vector  
bundle on $X$, which is canonically trivial on the boundary $\pa X$.  
The general definition of the index of the family $T$ in \cite{AS4} is  
that of the difference of the kernel bundle $\ker (V)$ and the trivial  
bundle $X\times \CC^N$. Since the bundle $\ker (V)$ is canonically  
trivial on the boundary of $X$, we obtain an element  
\begin{equation}\label{eq.def.deg}  
	\deg(T) \in K^0(X,\pa X).  
\end{equation}  
The degree is invariant with respect to homotopies $T_t$ of families  
of operators on $X$ that are {\em invertible} on $\pa X$ throughout  
the homotopy. We shall use the degree in Theorem \ref{Theorem.degree}  
for $X=B_R$, a large closed ball in $\GR^* \simeq \lgg^*$, or for $X$  
being the radial compactification of $\lgg^*$. If the boundary of $X$  
is empty, this construction goes back to Atiyah and Singer and then  
$\deg(T)$ is simply the family index of $T$. 
This definition of $\deg(T)$ is due to Melrose.  We note that when  
$\pa X \not = \emptyset$, the degree is not a local quantity in $T$,  
in the sense that it depends on more than just the principal symbol.  
  
Assume now that the family $T$ above consists of order zero operators  
and $T(x)$ is a multiplication operator for each $x \in \pa X$. We  
want to compute the Chern character of $\deg(T)$ using the  
Atiyah-Singer family index formula \cite{AS1,AS4}.  To introduce the  
main ingredients of the index formula, denote by $S^*_{vert}M$ the set  
of unit vectors in the dual of the vertical tangent bundle $T_{vert}M$  
to the fibers of $M \to X$.  Because the family $T$ is elliptic, the  
principal symbols define an invertible matrix of functions  
$$  
	a = \sigma_0(T) \in \CI(S^*_{vert}M;\Hom(E_0,E_1)).  
$$   
Since the operators $T(x)$ are multiplication operators, we can then  
extend $a$ to an invertible endomorphism on  
$$
	S_M:=S^*_{vert}M \cup B^*,
$$
with $B^*$ denoting the set of vertical  
cotangent vectors of length $\le 1$ above $\pa X$, as in  
\cite{AB}.   
The constructions of \cite{AS1,AS4} are in terms of    
$[a']\in K^0(T^*_{vert}M, T^*_{vert} M\vert_{\pa M})$ obtained  
by applying the clutching (or difference) construction to  
$a$. Explicitly, $[a']$ is represented by $(E_0,E_1,a_1)$ (where $a_1$ is  
a smooth function that coincides with $a$ outside a neighborhood of the  
zero section). It defines an element in  
$$  
	K^0(T^*_{vert}M, T^*_{vert} M\vert_{\pa M}) = K^0(T^*_{vert}M  
	\smallsetminus T^*_{vert} M\vert_{\pa M}),  
$$   
because $a$ defines an endomorphism of the trivial bundle $\CC^N$  
which is invertible outside a compact set. (Recall from  \cite{AtiyahKTHRY}
that we can represent the $K$-theory groups of a space $X$ as equivalence
classes of triples $(E_0,E_1,a)$, where $E_0$ and $E_1$ are vector bundles
on $X$ and $a$ is a vector bundle morphism $E_0 \to E_1$ which is an 
isomorphism outside a compact set).   
  
When $E_0 \cong E_1$, we can assume that  
$E_0 = E_1 = M \times \CC^N$ are trivial of rank $N$, and we have  
that $a$ is an invertible matrix valued function, which hence  
defines an element $[a] \in K^1(S_M)$.    
Let $B_1$ be the set of vectors of norm at most  
$1$ in $T^*M$. After the identification (up to homeomorphism) of $B_1  
\smallsetminus S_M$, the interior of $S_M$, with the difference  
$T^*_{vert}M \smallsetminus T^*_{vert} M\vert_{\pa M}$, we have  
$$  
	[a'] = \pa [a].  
$$  
  
We denote by $\pi_* : H^*(S_M) \to H^{* - 2n + 1}(X, \pa X)$ the  
integration along the fibers, where $n$ is the dimension of the fibers  
of $M \to X$. Integration along the fibers in this case is the  
composition of  
$$  
	\pa : H^*(S_M) \to H^{* + 1} (B_1, S_M) \cong  
	H_c^{*+1}(T^*_{vert}M \smallsetminus T^*_{vert} M\vert_{\pa M})  
$$   
and   
$$  
	\tilde\pi_*: H_c^*(T^*_{vert}M \smallsetminus T^*_{vert}  
	M\vert_{\pa M}) \to H_c^{* -2n}(X \smallsetminus \pa X)  
$$  
obtained by integration along the fibers of the bundle $T^*_{vert}M  
\smallsetminus T^*_{vert} M\vert_{\pa M} \to X \smallsetminus \pa X$:  
$$  
	\pi_* = \tilde \pi_* \circ \pa.  
$$  
  
To state the following result, we also need $Ch: K_1(S_M) \to  
H^{odd}(S_M)$, the Chern character in $K$-Theory and $\mathcal T$, the  
Todd class of $(T^*_{vert}M) \otimes \CC$, the complexified vertical  
tangent bundle of the fibration $M \to X$, as in \cite{AS4}.  Using  
the notation introduced above, we have:

\begin{theorem}\label{Theorem.mult.b}\ Let $M \to X$ be a smooth  
fiber bundle whose fibers are smooth manifolds (without corners), and  
let $T$ be a family of order zero elliptic pseudodifferential  
operators acting along the fibers of $M \to X$. Assume $X$ is a  
manifold with boundary $\pa X$ such that the operators $T(x)$ are  
multiplication operators on $\pa X$, also let $[a']$ and $[a]$ be the  
classes defined above. Then  
\begin{equation*}  
	Ch(\deg(T)) = (-1)^n \tilde\pi_*\big (Ch[a'] {\mathcal  
	T}\big) \in H^{*}(X, \pa X),  
\end{equation*}  
If $E_0 \cong E_1$, then we also have $Ch(\deg(T))= (-1)^n \pi_*\big  
(Ch[a] {\mathcal T}\big)$.  
\end{theorem}

\begin{proof}\ For continuous families $T(x)$ that are multiplication   
operators on the boundary $\pa X$, the degree is a local quantity --  
it depends only on the principal symbol -- so we can follow word for  
word \cite{AS4} to prove that  
$$  
        Ch(\deg(T)) = (-1)^n \tilde\pi_*\big (Ch[a'] {\mathcal  
        T}\big) \in H^{*}(X, \pa X).  
$$  
  
When $E_0 \cong E_1$, using $Ch [a'] = Ch (\pa [a]) = \pa Ch [a]$, we get  
\begin{multline*}  
	Ch(\deg(T)) = (-1)^n \tilde\pi_* \big ( \pa Ch[a] {\mathcal  
      	T} \big )  = (-1)^n \tilde \pi_* \circ \pa (Ch[a]  
      	{\mathcal T}) \\  = (-1)^n \pi_*\big (Ch[a] {\mathcal  
      	T}\big) \in H^{*}(X, \pa X).  
\end{multline*}  
This completes the proof.  
\end{proof}

It is interesting to note that it is not possible in general to give a  
formula for $Ch(\deg(T))$ only in terms of its principal symbol,  
without further assumptions on $T$ (for example that $T$ consists of 
multiplication operators on the boundary, as in our theorem). 
A consequence is that, in general,  
the formula for $Ch(\deg(T))$ will involve some non-local  
invariants. It would be nevertheless useful to find such a formula.  
   
Returning to our considerations, we continue to assume that $Z:=Y  
/\GR$ is compact, and we fix a metric on $\GR$ (which, we recall, is a  
vector bundle in these considerations). If $A \in \PsS{\infty}{Y;E_0,E_1}$   
is elliptic (in the sense that its principal symbol is invertible outside  
the zero section), then the indicial operators $\hat{A}(\tau)$ are  
invertible for $|\tau|\ge R$, $\tau \in \GR^*$, and some large $R$.  
In particular, by restricting the family $\hat{A}$ to the ball  
$$   
        B_R := \{ |\tau| \le R \},   
$$  
we obtain a family of elliptic operators that are invertible on the  
boundary of $B_R$, and hence $\hat{A}$ defines an element   
\begin{equation} \label{eq.def.degree}  
        \deg_{\GR}(A): = \deg(\hat{A}) \in K^0(B_R,\pa B_R) \simeq  
        K^0(\lgg^*)  
\end{equation}   
in the $K$-group of the ball of radius $R$, relative to its boundary,  
as explained above, called also the {\em degree} of $A$.  
  
We want a formula for the Chern character of the degree of   
$A \in \PsS {\infty}{Y;E_0,E_1}$.  
Because $A$ is elliptic, its principal symbol  
defines a class $[a'] \in K^0((T^*_{vert}Y)/\GR)$. If $E_0$ and  
$E_1$ are isomorphic, then it also defines a class   
$[a] \in K^0((S^*_{vert}Y)/\GR)$.   
Denote by $n$ the dimension of the quotient $Z_b=Y_b /\GR_b$ (which is  
independent of $b$ because we assumed $B$ connected), and  
let $\pi: (S^*_{vert}Y)/\GR \to B$ be the projection and  
\begin{equation}  
\begin{gathered}  
	\tilde\pi_* : H^*((T^*_{vert}Y)/\GR) \to H_c^{* - 2n }(\lgg^*)  
	\\ \pi_* : H^*((S^*_{vert}Y)/\GR) \to H_c^{* - 2n + 1}(\lgg^*)  
\end{gathered}  
\end{equation}  
be the integration along the fibers in cohomology.  Then  
$$  
	Ch(\deg_\GR(A)) \in H_c^*(\lgg^*)   
	\simeq H^{* +n}(B, \mathcal O),  
$$  
and the following theorem gives a formula for this cohomology class in  
terms of the classes $[a']$ or $[a]$ defined above.  
  
We denote by ${\mathcal T}$ the Todd class of the vector bundle  
$(T_{vert}Y)/\GR \otimes \CC \to Y/\GR$. We assume $B$ to be compact.

\begin{theorem}\label{Theorem.degree}\   
If $A \in \PsS{\infty}{Y;E_0,E_1}$ is elliptic, then the Chern  
character of $\deg_\GR(A)$ is given by  
\begin{equation*}  
	Ch( \deg_\GR(A)) = (-1)^n \tilde\pi_* \big (Ch[a'] {\mathcal T}\big)   
	\in H_c^{*}(\lgg^*),   
\end{equation*}  
Moreover, $Ch( \deg_\GR(A)) = (-1)^n \pi_* \big (Ch[a] {\mathcal  
T}\big)$ if $E_0 \cong E_1$.  
\end{theorem}

{\em Observations.} It is almost always the case that $E_0 \cong E_1$.  
For example, it is easy to see that this must happen if the Euler  
characteristic of $(T^*_{vert}Y)/\GR$ vanishes. This assumption is satisfied if 
$\GR = B \times \RR^q$, $q > 0$, which is the case for parameter dependent
pseudodifferential operators, for example.  
  
Another observation is that if the set of elliptic elements in  
$\PsS{\infty}{Y;E_0,E_1}$ is not empty, then $Z : = Y/\GR$ is compact.  
\vspace*{2mm}

\begin{proof}\ We cannot use Theorem \ref{Theorem.mult.b} directly  
because our family $\hat{A}$ does not consist of multiplication  
operators on the boundary. Nevertheless, we can deform $\hat{A}$ to a  
family of operators that are multiplication operators at $\infty$, for  
suitable $A$. We now construct this deformation.  
  
Let $E$ be a vector bundle over $V$. We consider classical symbols  
$S^m_c(E)$ whose support projects onto a compact subset of $V$.  Let  
$$  
	A_Y := (T_{vert}Y)/\GR \cong T_{vert}Z \times_B \lgg.  
$$  
First we need to define a nice quantization map $q : S^{m}(A_Y^*) \to  
\Psm{m}Y$. To this end, we proceed as usual, using local coordinates,  
local quantization maps, and partitions of unity, but being careful to  
keep into account the extra structure afforded by our settings:\ the  
fibration over $B$ and the action of $\GR$. Here are the details of  
how this is done.  
  
Fix a cross section for $Y \to Z := Y/\GR$ and, using it, identify $Y$  
with $Z \times_B \GR$ as $\GR$-spaces. Denote by $p_0 : Z \to B$ the  
natural projection. We cover $B$ with open sets $U_{\alpha'}$ that are  
diffeomorphic to open balls in a Euclidean space such that $Z  
\vert_{U_{\alpha'}} \cong U_{\alpha'} \times F$ and $\GR  
\vert_{U_{\alpha'}} \cong U_{\alpha'} \times \RR^q$. We also cover $F$  
with open domains of coordinate charts $V_{\alpha''}$.  Then we let  
$W_{\alpha} = p_0^{-1}(U_{\alpha'}) \cap V_{\alpha''}$, with $\alpha =  
(\alpha', \alpha'')$. The natural coordinate maps on $W_{\alpha}$ then  
give rise to a quantization map  
\begin{equation}  
	q_\alpha: S_c^m(A^*_Y\vert_{W_\alpha}) \to  
	 \Psm{m}{W_\alpha \times_B \GR}, \quad   
	q_\alpha(a) = a(b,x,D_x,D_t),  
\end{equation}  
where we identify  
\begin{equation}  
	S_c^m(A^*_Y\vert_{W_\alpha}) = S_c^m(T^*_{vert}(W_\alpha) \times_B \GR) =  
	S_c^m(U_{\alpha'} \times T^*V_{\alpha''} \times {\RR^*}^q).  
\end{equation}  
Denote by $b \in U_{\alpha'}$, $(x,y) \in T^*V_{\alpha''} \cong  
V_{\alpha''} \times {\RR^*}^n$, and $\tau \in {\RR^*}^q$ the  
corresponding coordinate maps. Then $q_\alpha(a) = a(b,x,D_x,D_t)$  
acts on $\CIc(U_{\alpha'} \times T^*V_{\alpha''} \times \RR^q)$ as  
\begin{multline*}  
	a(b,x,D_x,D_t)u(b,x,t) \\ = (2\pi i)^{-n - q} \int_{\RR^{n +  
	q}} \left (\int_{\RR^{n + q}}e^{i(x - z)\cdot y + i (t -  
	s)\cdot \tau} a(b,x,y,\tau) u(b,z,s) ds dz \right ) dy d\tau.  
\end{multline*}   
Choose now a partition of unity $\phi_\alpha^2$ subordinated to  
$W_\alpha$ and let  
$$  
	q(a) = \sum_\alpha q_\alpha(\phi_\alpha a) \phi_\alpha.  
$$   
  
The main properties of $q$ are the following:

\begin{enumerate}  
\item\ if $a$ has order $m$, then $\sigma_m(q(a)) = a$, modulo symbols of  
lower order;  
\item\ there exist maps $q_b : S^m(T^* Z_b ) \to \Psi^m(Z_b)$ such  
that  
$$  
	\widehat{q(a)}(\tau) = q_b(a(\cdot,\tau)) = a(b,x,\tau + D_x),  
$$  
if $\tau \in \lgg_b^*$ and $x \in Z_b$;  
\item\ the maps $q_b$ define a quantization map   
$$  
	\tilde q : \Gamma( \lgg^*, S^m( T^*_{vert}Z)) \to  
	\Gamma(\lgg^*, \Psi^m(Z)),  
$$  
where we regard $\Psi^m(Z_b)$ as defining a bundle of algebras,  
$\Psi^m(Z)$, on $B$, first, and then on $\lgg^*$, by pull-back.  
Similarly, we regard $S^m(T^*_{vert}Z)$ as defining a bundle over $Z$,  
which we then pull back to a bundle on $\lgg^*$.  (See also the  
discussion related to Equations \eqref{eq.bundle.hatA} and  
\eqref{eq.section.hatA}.)  
\end{enumerate}

The deformation of our family is obtained as follows. Let $|\tau|$ and 
$|y|$ be the norms $\tau \in \lgg^*$ and $y \in T^*_{vert}Z$. Define 
then 
\begin{equation}  
	\phi_\lambda^2 (y,\tau) = 1 + \lambda |y|^2 + \lambda (1 + 
	\lambda |\tau|^2)^{-1} |y|^2 \quad \text{and }\; 
	\psi_\lambda^{-2} (y,\tau) = 1 + \lambda |\tau|^2, 
\end{equation}  
which were chosen to satisfy 
$$
	1 + \phi_\lambda^2 |\tau|^2 +  
	\psi_\lambda^2 |y|^2 = 1 + |\tau|^2 + |y|^2 + \lambda |\tau|^2|y|^2.
$$  
For any symbol $a \in S^0(A_Y^*)$, $A_Y^* = T^*_{vert}Z \times_B  
\lgg^*$, we let  
$$  
	a_{\lambda,\tau}(y) = a( \psi_\lambda y, \phi_\lambda \tau ),  
	\quad \lambda \in [0,1],\; \tau \in \lgg^*_b, \text{ and } y  
	\in T^*Z_b.  
$$  
We can define then $A_{\lambda}(\tau) := q_b(a_{\lambda,\tau})$, $\tau  
\in \lgg^*_b$, which is the same as saying that $A_\lambda =  
\tilde{q}(a_{\lambda,\tau})$, and these operators will satisfy the  
following properties:  
  
\begin{enumerate}   
\item\ For each fixed $\lambda$, the operators $A_{\lambda}(\tau)$  
define a section of $\Psi^m(Z)$ over $\lgg^*$ and these sections  
depend smoothly on $\lambda$ (in other words, $A_{\lambda}(\tau)$  
depends smoothly on both $\lambda$ and $\tau$, in any trivialization);  
\item\ $A_{0}(\tau) = \widehat{q(a)}(\tau)$, for all $\tau$;  
\item\ For each nonzero $\tau' \in \lgg^*$ and $\lambda >0$, the limit  
$\displaystyle{\lim_{t \to \infty}} A_{\lambda}(t\tau')$ exists and is  
a multiplication operator;  
\item If $a,b \in S^0(A_Y^*)$ are such that $ab = 1$, $a$ is  
homogeneous of order zero outside the unit ball, and if we define  
$A_{\lambda} := \tilde{q}(a_{\lambda,\tau})$ and $B_{\lambda,\tau} :=  
\tilde{q}(b_{\lambda,\tau})$, then there exists a constant $C >0$ such  
that  
$$  
	\| A_\lambda(\tau) B_\lambda(\tau) - 1 \| \le C(1 + |\tau|)^{-1}  
$$  
and similarly  
$$  
	\| B_\lambda(\tau) A_\lambda(\tau) - 1 \| \le C(1 + |\tau|)^{-1}  
$$  
For all $\tau$ and $\lambda$;  
\item All these estimates extend in an obvious way to matrix valued  
symbols.  
\end{enumerate}

These properties are proved as follows. We first recall that, for any 
vector bundle $E$, we can identify the space of classical symbols 
$S_c^0(E)$ with $\CIc(E_1)$, the space of compactly supported 
functions on $E_1$, the unit ball of $E$, by $E_1 \smallsetminus \pa 
E_1 \ni \xi \to (1 - \|\xi\|^2)^{-1}\xi \in E$. For any quantization 
map, the norm of the resulting operator will depend on finitely many 
derivatives. Because we can extend $a_{\lambda,\tau}$ to a smooth 
function on the radial compactification of $A_Y^*$, the first property 
follows. The second property is obvious.  The third property is 
obtained using the same argument and observing that, for $\lambda >0$, 
we can further extend our function to the radial compactification in 
$\lambda$ also. By investigating what this limit is along various 
rays, we obtain the third property. 
  
The fourth property is obtained using the following  
observation: there exist a constant $C>0$ and seminorms $\|\,\cdot\,\|_0 $ on  
$S_c^0(T^*\RR^n)$ and $\|\,\cdot\,\|_{-1}$ on $S^{-1}(T^*\RR^n)$  
such that, for any symbols $a,b \in S_c^0(T^*\RR^n)$,  
$$  
	\|(ab)(x,D_x) - a(x,D_x) b(x,D_x) \| \le \sum_j    
	C(\|a\|_0 \| \pa_{y_j} b\|_{-1}  
	+ \| \pa_{y_j} a\|_{-1}\|b\|_0 ),  
$$  
$\pa_{y_j}$ being all derivatives in the symbolic directions (whose  
coordinates are denoted by $y$).   
Finally, the fifth property is obvious.  
  
We now turn to the proof of the formula for the degree of $A$ stated  
in our theorem. We prove it by a sequence of successive reductions,  
using the facts established above. First, it is easy to see that  
$\deg_\GR(A)$ depends only on its principal symbol, and hence we can  
assume that $A$ has order zero and $A = q(a)$, where $a =  
\sigma_0(A)$.  
  
The above deformation can be used to prove our theorem as follows.  
Fix $R$ large enough, and restrict the families $A_{\lambda}$ to the  
closed ball of radius $R$ in $\lgg^*$. For $|\tau| = R$ large enough,  
all operators $A_{\lambda}(\tau)$, $\lambda \in [0,1]$ are invertible,  
so the degree $\deg (A_\lambda)$ of these families is defined and does  
not depend on $\lambda$ or $R$, provided that $R$ is large  
enough. Since $\deg_\GR(A) = \deg(A_0)$, by definition, it is enough  
to compute $\deg(A_\lambda)$, for any given $\lambda$.  Choose then  
$\lambda >0$ arbitrary, and let $R \to \infty$.  Then the family  
$A_\lambda$ extends to a continuous family on the radial  
compactification of $\lgg^*$, which consists of multiplication  
operators on the boundary. Moreover, the symbol class of $A_\lambda$  
is nothing by the extension of $a_{\lambda}(\tau)$ to the radial  
compactification in $\tau$ and $\lambda$ (which is a manifold with  
corners of codimension two).  
  
We can use then Theorem \ref{Theorem.mult.b} to conclude that   
$$  
	\deg(A_\lambda) = (-1)^{n} \tilde\pi_* ( Ch [a_\lambda']  
	{\mathcal T}) \in H^*_c(\lgg^*).  
$$  
But $a_\lambda$ is homotopic to $a$ through symbols that are  
invertible outside a fixed compact set, so $[a_\lambda']=[a']$.  We  
get  
$$  
	\deg(A) = (-1)^{n} \tilde\pi_* ( Ch [a']  
	{\mathcal T}) \in H^*_c(\lgg^*).  
$$  
  
To obtain the second form of our formula for $E_0 \cong E_1$,  
and thus finish the proof, we proceed as at the end of the proof of  
Theorem \ref{Theorem.mult.b}, using $Ch [a'] = \pa Ch [a]$.    
\end{proof}

To prove the following result, we shall use terminology from algebraic  
topology:\ if $I_k \subset A_k$ are two-sided ideal of some algebras  
$A_0$ and $A_1$ and $\phi:A_0\to A_1$ is an algebra morphism, we say  
that $\phi$ {\em induces a morphism of pairs} $\phi:(A_0,I_0) \to  
(A_1,I_1)$ if, by definition, $\phi(I_0) \subset I_1$.

\begin{theorem}\label{Theorem.inddeg}\ Let $\GR$ be a vector bundle 
and $A \in \PsS m{Y;\CC^N}$ be an elliptic element. Then  
$$   
        \ind_\GR(A) = \deg_{\GR}(A) \in K^0(\lgg^*).   
$$   
\end{theorem}

\begin{proof}\  By embedding $E$ in a trivial vector bundle,
we can reduce to the case $E = B \times \CC^N$.
Let $B_R = \{ |\tau| \le R\} \subset \GR^*$ be as above.  The  
algebra   
$$  
      {\mathfrak A}_R := \CI(B_R,\Ps {\infty}{Y_b})  
$$  
of $\CI$-families of pseudodifferential operators on $B_R$ acting on  
fibers of $Y \times_B B_R \to B_R$, contains as an ideal ${\mathfrak  
I}_R = \CI_0(B_R,\Ps {-\infty}{Y_b})$, the space of families of  
smoothing operators that vanish to infinite order at the boundary of  
$B_R$.  If $A$ is an elliptic family, as in the statement of the  
lemma, and if $R$ is large enough, then $\hat{A}$, the indicial family  
of $A$, defines by restriction an element of $M_N({\mathfrak A}_R)$  
that is invertible modulo $M_N({\mathfrak I}_R)$.  
  
Recall that the boundary map $\pa_1$ in algebraic $K$-theory  
associated to the ideal ${\mfk I}_R$ of the algebra ${\mathfrak A}_R$  
gives $ \pa_1[A]= \deg_{\GR}(A),$ by definition. Also, the boundary  
map $\pa_0$ in algebraic $K$-theory associated to the ideal  
$\PsS{-\infty}Y$ of the algebra $\PsS {\infty}{Y}$ gives  
$\pa_0[A]=\ind_\GR(A)$. We want to prove that $\pa_1[A]=\pa_0[A]$.  The  
desired equality will follow by a deformation argument, which involves  
constructing an algebra smoothly connecting the ideals ${\mfk I}_R$  
and $\PsS {-\infty}Y$.  
  
Consider inside $\CI([0,R^{-1}], \PsS{-\infty}Y)$ the subalgebra of  
families $T=(T_x)$ such that $\hat{T}_x(\tau) = 0$ for $|\tau| \ge  
x^{-1}$. (In other words, $T_x \in {\mfk I}_{x^{-1}}$, if $x \not =  
0$, and $T_{0}$ is arbitrary.)  Denote this subalgebra by ${\mfk  
I}_{R\infty}$. Also, let ${\mfk A}_{R\infty}$ be the set of families  
$A=(A_x)$, $x \in [0,R^{-1}]$, $A_x \in {\mfk A}_{x^{-1}}$, if $x \not  
=0$, $A_0 \in \PsS {\infty}Y$ arbitrary such that the families $AT :=  
(A_xT_x)$ and $TA =: (T_x A_x)$ are in ${\mfk I}_{R\infty}$, for all  
families $T=(T_x) \in {\mfk I}_{R\infty}$.  
  
It follows that ${\mfk I}_{R\infty}$ is a two-sided ideal in ${\mfk  
A}_{R\infty}$ and that the natural restrictions of operators to  
$x=R^{-1}$ and, respectively, to $x=0$, give rise to morphisms of  
pairs  
$$  
       e_{1}:({\mfk A}_{R\infty}, {\mfk I}_{R\infty}) \to   
       ({\mfk A}_{R}, {\mfk I}_{R}) , \quad \text{ and }  
$$  
$$  
       e_{0}: ({\mfk A}_{R\infty}, {\mfk I}_{R\infty}) \to (\PsS  
       {\infty}Y, \PsS {-\infty}Y).  
$$  
Moreover, the indicial family of the operator $A$ gives rise, by  
restriction to larger and larger balls $B_r$, to an invertible element  
in ${\mfk A}_{R\infty}$, also denoted by $A$. Let $\pa$ be the  
boundary map in algebraic $K$-theory associated to the pair $({\mfk  
A}_{R\infty},{\mfk I}_{R\infty})$. Then $(e_0)_* \pa[A]= \pa_0[A]$ and  
$(e_1)_* \pa[A]=\pa_1[A]$. Since $(e_0)_*: K_0({\mfk I}_{R\infty}) \to  
K_0(\PsS {-\infty}Y)$ and $(e_1)_*: K_0({\mfk I}_{R\infty}) \to  
K_0({\mfk I}_{R})$ are natural isomorphisms, our result follows.  
\end{proof}

We now drop the assumption above that $\GR$ consist of abelian Lie 
groups, assuming instead that $\GR$ consists of simply-connected 
solvable Lie groups, and want to compute the Chern character of the 
gauge-equivariant  index $\ind_\GR(A)$, for an elliptic family $A \in \Psm {m}Y$. 
One difficulty that we encounter is that the space on which the 
principal symbols are defined, that is $(S^*_{vert}Y)/\GR$, is not 
orientable in general. (Recall that $S^*_{vert}Y$ is the space of 
vectors of length one of $T_{vert}^*Y$, the dual of the vertical 
tangent bundle $T_{vert}Y$ to the fibers of $Y \to B$.) 
  
We denote by ${\mathcal T}$ the Todd class of the vector bundle  
$(T_{vert}Y)/\GR \otimes \CC \to Y/\GR$ and by $\pi_*$ the integration  
along the fibers of $(S^*_{vert}Y)/\GR \to B$, as above. We assume $B$   
to be compact.

\begin{theorem}\label{Theorem.Chern}\   
Let $\GR$ be a bundle of Lie groups whose fibers are simply-connected,  
solvable Lie groups. Let $A \in \Psm {m}{{Y,E}}$ be an elliptic,  
$\GR$-invariant family, and let $[\sigma_m(A)] \in K^1((S^*_{vert}Y)/\GR)$ be  
the class defined by the principal symbol $\sigma_m(A)$ of $A$. Then  
the Chern character of the gauge-equivariant  index of $A$ is given by  
$$   
        Ch(\ind_\GR(A)) = (-1)^n \pi_*\big (Ch[\sigma_m(A)]   
        {\mathcal T}\big) \in H_c^{*}(\lgg^*),  
$$  
where $n$ is the dimension of the fibers of $(S^*_{vert}Y)/\GR \to B$.  
\end{theorem}

\begin{proof}\      
Note first that we can deform the bundle of Lie groups $\GR$ to the  
bundle of {\em commutative} Lie groups $\lgg$ as before, using  
$\GR_{adb}$. Moreover, we can keep the principal symbol of $A$ constant  
along this deformation. This shows that we may assume $\GR$ to consist  
of commutative Lie groups, i.e. that $\GR$ is a vector  
bundle. The result then follows from Theorems \ref{Theorem.degree} and  
\ref{Theorem.inddeg}.  
\end{proof}

{\em Observations.} We can extend the above theorems in several ways.  
First, we can drop the assumption that $Z= Y/\GR$ be compact, but then  
we need to consider bounded, elliptic elements $A \in \Hom(E,E) +  
\Psm 0{Y;E}$ (or, if $\GR$ is a vector bundle, then 
$A \in \Hom(E,E) + \Psm 0{Y;E}$).  Also, in the last two theorems, we 
can allow operators acting between sections of different vector 
bundles. This will require to slightly modify the proof of Theorem 
\ref{Theorem.inddeg}, either by using a smooth version of bivariant 
$K$-theory \cite{NistorKthry}, or by using the usual bivariant 
$K$-theory after we have taken the norm closures of the various ideals 
$\mfk I$ decorated with various indices. We can also further integrate 
along the fibers of $\lgg^* \to B$ to obtain a cohomological formula 
with values in $H_c^*(B;\mathcal O)$, the cohomology with local 
coefficients in the orientation sheaf of $\lgg^* 
\to B$.

\section{Index theory on a simplex\label{S.Index.Simplex}}

In this section, we discuss an application of our main index
theorem to the question of formulating
(\ie the existence of) Fredholm boundary conditions for 
$b$-pseudodifferential operators on the simplex
$$
        \Delta_n=\{(x_0,x_1,\ldots,x_n), x_i \ge 0, \sum x_i=1\}
$$ 
and on manifolds with corners of the form $\Delta_n \times Y$, 
where $Y$ is a smooth, compact
manifold without corners.  We obtain complete results for $\dim Y >
0$.  The answer to this existence question is ``local,'' that is, it
can be given in terms of the principal symbol. The problem of
computing the index of the resulting Fredholm operators, if any, is a
``non-local'' problem. Some results of this section were obtained in
a conversation with R. Melrose, whose insights I am happy to acknowledge.

In addition to the results of the previous sections, we shall also use
computations from \cite{Melrose-Nistor1}. All the definitions not
included in this section can be found in that paper.

We formalize the above question as the following natural problem:\\[2mm]
{\bf Problem $({\mathcal F}T)$.}\ {\em 
Suppose an elliptic operator $T \in \Psi_b^m(\Delta_n\times Y;E)$ is
given. For what $T$ can we find a perturbation $T+R$, with 
$R \in \Psi_b^{-\infty}(\Delta_n)$, such that
$$
	T + R:H^s(\Delta_n\times Y;E)
	\to H^{s-m}(\Delta_n\times Y;E)
$$ 
is Fredholm?}\\[2mm]
We can reduce the general case to the case $m = 0$ by replacing $T$
with $D^{-m}T$, where $D$ is an elliptic, strictly positive operator
in $\Psi_b^1(\Delta_n \times Y;E)$. (The resulting operator $D^{-m} T$
is then in the norm completion of $\Psi_b^0(\Delta_n \times
Y;\End(E))$, but this makes no difference to us.) For $m=0$, we then
obtain an operator acting on the Hilbert space $L^2(\Delta_n)$.

The answer to the above problem, Problem $({\mathcal F}T)$, is
independent on the metric chosen on $\Delta_n \times Y$, although the
choice of $R$ may depend on the metric.

Let $\overline \Psi_b^0(\Delta_n \times Y)$ be the norm closure of
$\Psi_b^0(\Delta_n \times Y)$. By the results of
\cite{Melrose-Nistor1}, the algebra $\overline \Psi_b^0(\Delta_n
\times Y)$ contains the ideal of compact operators ${\mathcal K}$ on
$L^2(\Delta_n \times Y)$.  (See also \cite{LMN}.)  Denote by $Q$ be
the algebra of ``joint symbols'' of $\Psi_b^0(\Delta_n \times Y)$:
$$
        Q:=\overline \Psi_b^m(\Delta_n \times Y)/\mathcal K,
$$
The principal symbol then extends to a surjective morphism
$\sigma_0 : Q \to C(\bS^*\Delta_n \times Y)$ with kernel $J$. 

A first, tautological observation is that Problem $({\mathcal F}(T))$
has a solution if, and only if, the invertible element $\sigma_0(T)
\in C(\bS^*\Delta_n \times Y;\End(E))$ has a lifting to an invertible
element in $Q$ (the lifting is $T + R$). Since there is a necessary
condition for this to happen in terms of $K$-theory, and the
$K$-groups involved do not depend on $E$, up to isomorphism, we obtain
the following necessary condition for our problem to have a solution:
Problem $({\mathcal F}(T))$ has a solution only if the class
$[\sigma_0(T)] \in K_1(C(\bS^*\Delta_n \times Y))$ is in the image of
the morphism $K_1(Q) \to K_1(C(\bS^*\Delta_n \times Y))$.

Now the standard six term $K$-theory exact sequence associated to the
short exact sequence of $C_r^*$-algebras 
$$
     0 \to J \to Q \to C(\bS^*\Delta_n \times Y) \to 0
$$
tells us that the class $[\sigma_0(T)] \in K_1(C(\bS^*\Delta_n \times
Y))$ is in the image of the morphism $K_1(Q) \to K_1(C(\bS^*\Delta_n
\times Y))$ if, and only if, it is in the kernel of the boundary map
$\pa:K_1(C(\bS^*(\Delta_n \times Y))) \to K_0(J)$ associated to the
above exact sequence of algebras.

\begin{lemma}\label{Lemma.cond}
If Problem (${\mathcal F}_T$ has a solution, then $\pa [\sigma_m(T)] =
0 \in K_0(J)$. This condition is also sufficient if $\dim Y >0$.
\end{lemma}

\begin{proof}\ The first part was already proved.
The second part follows using the same reasoning as in the
proof of from Theorem \ref{Theorem.obst}.   
\end{proof}

Denote by $I=\overline \Psi_b^{-1}(\Delta_n \times Y)$ the norm
closure of $\Psi_b^{-1}(\Delta_n \times Y)$ (the closure is in the
norm topology of bounded operators on $L^2(\Delta_n\times Y)$). Then
$J=I / \mathcal K$.  The results of \cite{Melrose-Nistor1} give a
spectral sequence converging to the $K$--theory groups of $I$ (and,
with some obvious changes, a spectral sequence converging to the
$K$--theory of $J$). The $E^1$-term of this spectral sequence is
independent of $Y$ and is (dual to) the combinatorial simplicial
complex of $\Delta_n$, more precisely, it is a direct sum of complexes
isomorphic to
$$
	0 \lra \ZZ^n \lra \ZZ^{\frac{n(n-1)}{2}} 
	\lra \ldots \ZZ^n \lra \ZZ \lra 0.
$$
The direct sum commes from the fact that $K$-theory is indexed by
$\ZZ/2\ZZ$ and not by $\ZZ$. The resulting spectral sequence, which is
indexed by $\ZZ \times \ZZ$, will be periodic of period $2$ in each
variable.

It follows that the spectral sequence converging to $K_n(I)$
degenerates at $E^2$, which shows that $K_n(I)\simeq \ZZ$ and
$K_{n-1}(I)=0$ (the same $n$ as in $\Delta_n$). To obtain the
$K$--theory of $J$ one proceeds similarly, the only difference being
that one drops the last copy of $\ZZ$ in the above complex. Finally,
this gives $K_0(J) \simeq K_0(I)$ and $K_1(J) \simeq K_1(I) \oplus
\ZZ$.

Denote by 
$$
	\In_k: \Psi_b^\infty(\Delta_n \times Y) \to \PsS{\infty}Y,
$$
$k=0,\ldots,n,$ the indicial maps corresponding to the $(n+1)$-corners of
$\Delta_n \times Y$.  If $B$ is reduced to a point and $\GR =\RR^n$,
we shall denote $\deg_{\GR}=\deg_n$. We are now ready to formulate and
prove the main result of this section.

\begin{theorem}\ Let $T \in \Psi_b^{m}(\Delta_n \times Y)$, $m>0$,
be an elliptic operator. Assume $\dim Y > 0$. Then, for $n$ odd,
Problem $({\mathcal F}(T))$ always has a solution.  For $n$ even,
$\deg_n(\In_k(T))$ is independent of $k$ and Problem $({\mathcal
F}(T))$ has a solution if, and only if, $\deg_n(\In_k(T))=0$.
\end{theorem}

\begin{proof}\ Using the results of \cite{Melrose-Nistor1},
we first observe that, in
order for $T + R$ to be Fredholm, it is necessary that $\In_k(T+R)$ be
invertible for all $k$, so
$$
        \deg_n(\In_k(T))=\deg_n(\In_k(T+R))=0, \ \text{ for all }k.
$$ 
The vanishing of all $\deg_n(\In_k(T))$ is hence necessary for Problem
$({\mathcal F}(T))$ to have a solution.  This condition is
automatically satisfied for $n$ odd.

To prove that the vanishing of any of $\deg_n(\In_k(T))$ is also
enough for Problem $({\mathcal F}(T))$ to have a solution, we shall
show that all degrees $\deg_n(\In_k(T))$ are equal, $k = 1, \ldots ,
n$, and that $\pa[\sigma_0(T)] = 0$ if, and only if,
$\deg_n(\In_k(T))=0$. The result will follow then from Lemma
\ref{Lemma.cond}.

As discussed above, it is enough to consider the case $m=0$. Now
$$
        K_0(J)=E^2_{0n}= \ker (\ZZ^n \to \ZZ^{n(n-1)/2}), 
$$ 
so $K_0(J)=\{(p,p,\ldots,p)\} \subset \ZZ^n$. Because
$$
        \pa[\sigma_0(T)]=(\deg_n(\In_0(T)),\deg_n(\In_1(T)),\ldots,
        \deg_n(\In_n(T))) \in E_{0n}^0,
$$ 
we obtain that the degrees $\deg_n(\In_k(T)))$ are all equal and that
the vanishing of $\pa[\sigma_0(T)]$ is equivalent to the vanishing of
any of the degrees $\deg_n(\In_k(T)))$. As observed above, an
application of Lemma \ref{Lemma.cond} is now enough to complete the
proof.
\end{proof}

\providecommand{\bysame}{\leavevmode\hbox to3em{\hrulefill}\thinspace}

\end{document}